\documentclass{amsart}
\newtheorem{theorem}{Theorem}[section]
\newtheorem{lemma}[theorem]{Lemma}
\newtheorem{cor}[theorem]{Corollary}
\theoremstyle{definition}
\newtheorem{defn}[theorem]{definition}
\newtheorem{ex}[theorem]{Example}

\newtheorem{assum}[theorem]{Assumption}

\theoremstyle{remark}

\numberwithin{equation}{section}

\usepackage[colorlinks]{hyperref}

\begin{document}
	\title{Additivity of multiplicative (generalized) skew semi-derivations on rings}
	\author{Sk Aziz}
    \address{Department of Mathematics, Indian Institute of Technology Patna, Patna-801106}
    \curraddr{}
    \email{E-mail: aziz$\textunderscore$2021ma22@iitp.ac.in}
    \thanks{}

   \author{Arindam Ghosh}
   \address{Department of Mathematics, Government Polytechnic Kishanganj, Kishanganj-855116}
   \curraddr{}
   \email{E-mail: arindam.rkmrc@gmail.com}
   \thanks{}

   \author{Om Prakash$^{\star}$}
   \address{Department of Mathematics, Indian Institute of Technology Patna, Patna-801106}
   \curraddr{}
   \email{om@iitp.ac.in}
   \thanks{* Corresponding author}
	
	\subjclass{16W10, 16S50}
	
	\keywords{Additivity; Skew derivation; Semi-derivation; (Generalized) Skew Semi-derivation; Idempotent}
	
	\dedicatory{}
	
	\maketitle
\begin{abstract}
	In this paper, we introduce a new class of derivations that generalizes skew derivations and semi-derivations, and we call it ``skew semi-derivation". Furthermore, we present a study of the conditions under which this type of multiplicative derivation becomes additive.

	\end{abstract}
	
\section{Introduction}
Let $R$ be a ring. A map $f: R \rightarrow R$ is said to be additive if $f(x+y)=f(x)+f(y)$ for all $x,y \in R$.
In 1937, Jacobson \cite{jacobson1937abstract} formally introduced the notion of a derivation on an algebra and proved that any derivation of a semisimple algebra over a ring could be extended to the inner derivation on the underlying ring. An additive map $h: R \rightarrow R$ is said to be a derivation if $h(xy)=h(x)y+xh(y)$ for all $x, y \in R$. More details on derivations over rings and algebras can be found in \cite{posner1957derivations, weisfeld1960derivations, kadison1967derivations, vukman1996derivations}.

In 1985, Leroy \cite{leroy1985derivees} introduced the notion of skew derivation over a ring and further studied it over skew fields and prime rings \cite{kharchenko1992skew}. An additive map $\delta: R \rightarrow R$ is said to be a skew derivation if $\delta(xy)=\delta(x)y+\alpha(x)\delta(y)$ for all $x, y \in R$, where $\alpha$ is an automorphism on $R$. Interested readers can see more results on skew derivations in \cite{chuang2000identities, chou2009engel}.

In 1983, Bergen \cite{bergan1983derivations} introduced the concept of semi-derivations as an extension of derivations of a ring $R$. Recall that an additive map $D: R \rightarrow R$ is said to be a semi-derivation if $D(xy)=D(x)g(y)+xD(y)=D(x)y+g(x)D(y)$ and $D(g(x))=g(D(x))$ where $x, y \in R$ and $g$ is an arbitrary map on $R$. We refer \cite{chang1984semi, brevsar1990semiderivations} to the readers for the basic results on semi-derivations over rings.
\section{Preliminaries}

These above notions raise the question of whether a unified definition of semi-derivation and skew derivation exists. Towards this, we have an affirmative answer and introduce two new definitions below that generalize semi-derivations and skew derivations.

\begin{defn}
	An additive map $d: R \rightarrow R$ is said to be a skew semi-derivation if
		$$
	\begin{aligned}
		&d(x y) =d(x) g(y)+\alpha(x) d(y) =d(x) \alpha(y)+g(x) d(y)\\
		 &\text{and}~ d(g(x))=g(d(x))
		\end{aligned}
	$$
	 where $ x, y \in R$, $\alpha$ is an automorphism and $g$ is a map on $R$.
\end{defn}
An additive map $\Delta$ is said to be a generalized derivation if $\Delta(xy)=\Delta(x)y+xh(y)$ for all $x,y \in R$ and some derivations $h$ on $R$. Now, we introduce the notion of generalized skew semi-derivation by combining the concepts of generalized derivation and skew semi-derivation as follows:

		\begin{defn}
			 An additive map $f: R \rightarrow R$ is said to be a generalized skew semi-derivation if
	$$
	\begin{aligned}
		&f(x y) =f(x) g(y)+\alpha(x) d(y) =d(x) \alpha(y)+g(x) f(y) \\
		& \text { and } f(g(x)) =g(f(x))
	\end{aligned}
	$$
		
for all $ x,y \in R$, where $d$ is a skew semi-derivation, $\alpha$ is an automorphism and for some map $g$ on $R$.
\end{defn}
	
Note that semi-derivations and skew derivations are obvious examples of skew semi-derivations. However, the converse is not true.
Towards this, we have an example of a skew semi-derivation that is neither a semi-derivation nor a skew derivation with the same associated map.
Let $e_{ij}$ be the matrix whose $(i,j)$-th entry is $1$, and all other entries are $0$. We consider a full matrix ring $M_n(\mathbb{R})$ over the field of real numbers $\mathbb{R}$ where $n\geq 2$ is an integer.
		\begin{ex}
We define three maps $d$, $\alpha$, and $g$ on $M_n(\mathbb{R})$ are as follows:\\
\begin{itemize}
\item The map $d: M_n(\mathbb{R})\rightarrow M_n(\mathbb{R})$ is defined by	
	\begin{align*}
d  \left( \left( \begin{array}{ccc}
a_{11} & a_{12} \\
 a_{21} & a_{22} \end{array}  \right)\right)= \left( \begin{array}{ccc}
2a_{22} & 2a_{21} \\
2a_{12} & 2a_{11} \end{array}  \right),
\end{align*}
	where $a_{ij}\in \mathbb{R}$.
	
\item The map $\alpha: M_n(\mathbb{R})\rightarrow M_n(\mathbb{R})$ is defined as  	
\begin{align*}
\alpha  \left(\left( \begin{array}{ccc}
a_{11} & a_{12} \\
 a_{21} & a_{22} \end{array}  \right)\right)= \left( \begin{array}{ccc}
a_{22} & a_{21} \\
a_{12} & a_{11} \end{array}  \right).
\end{align*}

\item The map $g: M_n(\mathbb{R})\rightarrow M_n(\mathbb{R})$ is defined by
\begin{align*}
g  \left(\left( \begin{array}{ccc}
a_{11} & a_{12} \\
 a_{21} & a_{22} \end{array}  \right)\right)= \left( \begin{array}{ccc}
0 & 0 \\
0 & 0 \end{array}  \right).
\end{align*}

\end{itemize}
Then $\alpha$ is an automorphism, and $g$ is a simple map. Also, an easy computation shows that $d$ is a skew semi-derivation with the associated map $g$ and automorphism $\alpha$ on $M_n(\mathbb{R})$. But taking $x=e_{11}$ and $y=e_{12}$, we have
\begin{align*}
d(xy)=2e_{21}\neq 0=d(x)g(y)+xd(y).
\end{align*}
Hence, $d$ is not a semi-derivation with the associated map $g$. Also, taking $x=e_{21}$ and $y=e_{22},$ we have
\begin{align*}
d(xy)=0\neq 2e_{12}=d(x)y+\alpha(x)d(y).
\end{align*}
Therefore, $d$ is not a skew derivation with the associated automorphism $\alpha$.
	\end{ex}

 It is noted that every skew semi-derivation is a generalized skew semi-derivation but the converse is not true as shown in the following example.
\begin{ex}
Again, we consider four maps as follows:
\begin{itemize}
\item A map $f : M_n(\mathbb{R}) \rightarrow M_n(\mathbb{R})$ is defined by
\begin{align*}
f  \left (\left( \begin{array}{ccc}
a_{11} & a_{12} \\
 a_{21} & a_{22} \end{array}  \right)\right)= \left( \begin{array}{ccc}
3a_{22} & 3a_{21} \\
3a_{12} & 3a_{11} \end{array}  \right);
\end{align*}
\item $g : M_n(\mathbb{R}) \rightarrow M_n(\mathbb{R})$ is defined as
\begin{align*}
g  \left(\left( \begin{array}{ccc}
a_{11} & a_{12} \\
 a_{21} & a_{22} \end{array} \right) \right)= \left( \begin{array}{ccc}
a_{22} & a_{21} \\
a_{12} & a_{11} \end{array}  \right);
\end{align*}
\item $d : M_n(\mathbb{R}) \rightarrow M_n(\mathbb{R})$ is defined as
\begin{align*}
d  \left(\left( \begin{array}{ccc}
a_{11} & a_{12} \\
 a_{21} & a_{22} \end{array}  \right)\right)= \left( \begin{array}{ccc}
0 & 0 \\
0 & 0 \end{array}  \right);
\end{align*}
and \item $\alpha : M_n(\mathbb{R}) \rightarrow M_n(\mathbb{R})$ is defined by
\begin{align*}
\alpha  \left(\left( \begin{array}{ccc}
a_{11} & a_{12} \\
 a_{21} & a_{22} \end{array}  \right)\right)= \left( \begin{array}{ccc}
a_{11} & -a_{12} \\
-a_{21} & a_{22} \end{array}  \right),
\end{align*}
where $a_{ij} \in \mathbb{R}$.
\end{itemize}

Then $\alpha$ is an automorphism, $g$ is a map, and by easy computation, $d$ is a skew semi-derivation with the associated map $g$ and automorphism $\alpha$. Also, it can easily be seen that $f$ is a generalized skew semi-derivation with the associated map $g$, automorphism $\alpha$, and skew semi-derivation $d$. However, taking $x=e_{21}$ and $y=e_{22}$, we have
\begin{align*}
f(x)\alpha(y)+g(x)f(y) = 3e_{12} \neq 0 = f(xy).
\end{align*}
Therefore, $f$ is not a skew semi-derivation with the associated map $g$ and automorphism $\alpha$.
\end{ex}	

Now, let us remove the additivity condition from the above definitions. In that case, we call these maps multiplicative derivation, multiplicative skew derivation, multiplicative semi-derivation, multiplicative skew semi-derivation, multiplicative generalized derivation, and multiplicative generalized skew semi-derivation, respectively. It is an exciting problem to show under what conditions these multiplicative maps become additive. The first result was seen due to Martindale \cite{martindale1969multiplicative} in 1969, where he proved that every multiplicative bijective map on a ring $R$ with a non-trivial idempotent is additive. More related results can be found in \cite{aziz2023some, aziz2023additivity, rickart1948one, johnson1958rings, daif1991multiplicative, wang2009additivity, wang2011additivity,ferreira2014multiplicative, el2014multiplicativity, yadav2017additivity}. In 2022, Siddeeque and Khan \cite{siddeeque2022multiplicative} proved that a multiplicative semi-derivation on a ring $R$ satisfying certain conditions becomes additive. These results motivate us to introduce skew semi-derivation and generalized skew semi-derivation and study these maps over rings. In this paper, we are trying to find an answer to the following question:
``Under what conditions do multiplicative skew semi-derivation and multiplicative generalized skew semi-derivation become additive?"\\

Throughout this paper, let $R$ be a ring with a non-trivial idempotent element $e$. Let $e_1=e$ and $e_2=1-e_1$ where $1$ is the identity element of $R$. Then $R=R_{11} \oplus R_{12} \oplus R_{21} \oplus R_{22}$ where $R_{ij}=e_i R e_j$ for $i, j\in \{1,2\}$. 	

\begin{assum}
\label{asm1.1}
Let $R$ be a ring with the zero element $0$ and $g:R\rightarrow R$ be a map with $g(0)=0$ and for $a_{ij}, b_{ij} \in R_{ij})$, $i, j, k=1,2$ satisfies the following conditions: \\

\noindent	(i) If $\ (a_{k 1}+b_{k 2}) g\ (x_{i j})=0$ and $i \leq j$, for all $x_{i j} \in R_{i j}$, then
$$
\begin{cases}a_{k 1}=0 & \text { if } i=1, \\ b_{k 2}=0 & \text { if } i=2.\end{cases}
$$

\noindent	(ii) If $g\ (x_{i i})\ (a_{1 j}+b_{2 j})=0$, for all $x_{i i} \in R_{i i}$, then
$$
\begin{cases}a_{1 j}=0 & \text { if } i=1, \\ b_{2 j}=0 & \text { if } i=2.\end{cases}
$$
(iii) $g (a_{11}+b_{12})=g (a_{11})+g (b_{12})$.\\
(iv) $g (a_{11}+b_{21})=g(a_{11})+g (b_{21}).$
\end{assum}

The fascinating part of this work is introducing a new class of derivations that extends the concepts of skew derivations and semi-derivations, called skew semi-derivations. Further, a more general concept is introduced, called generalized skew semi-derivations, and suitable examples are provided to support the presented notions.

When studying these new concepts, one natural question arises: do multiplicative skew semi-derivations and generalized skew semi-derivations become additive under certain conditions? We answer this question in the affirmative and provide proof in Sections 2 and 3.
Moreover, Section 2 proves that any multiplicative skew semi-derivation over a ring $R$ with a non-trivial idempotent $e$ is additive under certain assumptions on the underlying ring and maps. Section 3 extends this result to generalized skew semi-derivations over the same ring with additional assumptions.
Our results extend the theory of derivations over rings and provide new insights into the structure of skew derivations and semi-derivations.

\section{Skew Semi-Derivation}
\begin{theorem}
\label{thm2.1}
Let $R$ be a ring and $g: R \rightarrow R$ be a map that satisfies Assumption \ref{asm1.1}.
If $ d: R \rightarrow R$ is a multiplicative skew semi-derivation with the associated map $g: R \rightarrow R$ and an automorphism $\alpha: R \rightarrow R$, then $d$ is additive on $R$.
	\end{theorem}
	
	Before proving Theorem \ref{thm2.1}, we prove several lemmas.	
	
	\begin{lemma}
		\label{lem2.2}
		$d(0)=0$.
	\end{lemma}	
	\begin{proof}
		$d(0)=d(0\cdot 0)=d(0)g(0)+\alpha(0)d(0)=0$ (By Assumption \ref{asm1.1}, $g(0)=0$ and as $\alpha$ is an automorphism, $\alpha(0)=0$).
	\end{proof}
	
	\begin{lemma}
		\label{lem2.3}
		Let $a_{ij} \in R_{ij}$. Then
		$$
		\begin{aligned}
			(i)~ d (a_{11}+a_{12} )=d (a_{11} )+ d(a_{12} );\\
			(ii)~ d(a_{11}+a_{21} )= d(a_{11} )+d (a_{21} );\\
			(iii)~ d(a_{22}+a_{12} )=d (a_{22} )+d (a_{12} );\\
			(iv)~ d (a_{22}+a_{21} )=d (a_{22} )+d (a_{21}); \\
			(v)~ d(a_{11}+a_{22} )=d (a_{11} )+d(a_{22} ).
		\end{aligned}
		$$
	\end{lemma}

	\begin{proof} Let $x_{11} \in R_{11}$. Then
		\begin{equation}
			\label{eq2.1}
			d ( (a_{11}+a_{12} ) x_{11} )=d (a_{11}+a_{12} ) g (x_{11} )+\alpha (a_{11}+a_{12} ) d (x_{11} ).
		\end{equation}
		
		On the other hand,
		\begin{equation}
			\label{eq2.2}
			\begin{aligned}
				&d ( (a_{11}+a_{12} ) x_{11} )\\
				 &= d (a_{11} x_{11} )\\
				 &= d (a_{11} x_{11} )+d (a_{12} x_{11} ) \\
				&=d (a_{11} ) g (x_{11} )+\alpha (a_{11} ) d (x_{11} )+d (a_{12} ) g (x_{11} )  +\alpha (a_{12} ) d (x_{11} )
			\end{aligned}
		\end{equation}
		
		Comparing \eqref{eq2.1} and \eqref{eq2.2}, we have
		\begin{equation}
			\label{eq2.3}
			\begin{aligned}
				& { [d (a_{11}+a_{12} )-d (a_{11} )-d (a_{12} ) ] } g (x_{11} )=0 ~\text{which implies that,}\\ \\
				& ( [d (a_{11}+a_{12} )-d (a_{11} )-d (a_{12} ) ]_{11} \\&+  { [d (a_{11}+a_{12} )-d (a_{11} ) } -d (a_{12} ) ]_{12} ) g (x_{11} )=0\\
			\text{and}~	& ( [d (a_{11}+a_{12} )-d (a_{11} )-d (a_{12} ) ]_{21} \\&+  { [d (a_{11}+a_{12} )-d (a_{11} )} -d (a_{12} ) ]_{22} ) g (x_{11} )=0.
			\end{aligned}
		\end{equation}
		
		Similarly, for $x_{22} \in R_{22}$,
		\begin{equation}
			\label{eq2.4}
			\begin{aligned}
				&d (a_{11}+a_{12}) g (x_{22} )+\alpha (a_{11}+a_{12} ) d (x_{22}) \\
				&=d ( (a_{11}+a_{12} ) x_{22} )\\
				&=d (a_{12} x_{22} ) \\
				&=d (a_{11} x_{22} )+d (a_{12} x_{22} ) \\
				&=d (a_{11} ) g (x_{22} )+\alpha (a_{11} ) d (x_{22} )+d (a_{12} ) g (x_{22} )+\alpha (a_{12} ) d (x_{22} ).
			\end{aligned}
		\end{equation}	
		
		Comparing both sides of \eqref{eq2.4}, we have
		\begin{equation}
			\label{eq2.5}
			\begin{aligned}
				&\ [d (a_{11}+a_{12} )-d (a_{11} )-d (a_{12} ) ] g (x_{22} )=0~\text{which implies that,}\\ \\
				& ( [d (a_{11}+a_{12} )-d (a_{11} )-d (a_{12} ) ]_{11}\\& +[d (a_{11}+a_{12} )-d (a_{11} )-d (a_{12} ) ]_{12} ) g (x_{22} )=0\\
			\text{and}~	& ( [d (a_{11}+a_{12} )-d (a_{11} )-d (a_{12} ) ]_{21} \\&+  { [d (a_{11}+a_{12} )-d (a_{11} ) }-d (a_{12} ) ]_{22} ) g(x_{22} )=0 .
			\end{aligned}
		\end{equation}
		
		By the Assumption \ref{asm1.1} (i), \eqref{eq2.3} and \eqref{eq2.5}, we get
		$$
		\begin{aligned}
			&{ [d (a_{11}+a_{12} )-d (a_{11} )-d (a_{12} ) ]_{k 1} } =0 \\
			\&~ & { [d (a_{11}+a_{12} )-d (a_{11} )-d (a_{12} ) ]_{k 2 }}=0
		\end{aligned}
		$$
		for $k=1,2$. This completes the proof of Lemma \ref{lem2.3} (i). \\
		
		To proof of (ii) of Lemma \ref{lem2.3}, take $x_{11} \in R_{11}$, then
		\begin{equation}
			\label{eq2.6}
			\begin{aligned}
				& d (x_{11} ) \alpha (a_{11}+a_{21} )+g (x_{11} ) d (a_{11}+a_{21} )=d (x_{11} (a_{11}+a_{21} ) ) \\
				&=d (x_{11} a_{11} ) \\
				&=d (x_{11} a_{11} )+d (x_{11} a_{21} ) \\
				&=d (x_{11} ) \alpha (a_{11} )+g (x_{11} ) d (a_{11} )+d (x_{11} ) \alpha (a_{21} ) +g (x_{11} ) d (a_{21} ).
			\end{aligned}
		\end{equation}	
		
		Comparing both sides of \eqref{eq2.6}, we get
		\begin{equation}
			\label{eq2.7}
			\begin{aligned}
				&g (x_{11} ) [d (a_{11}+a_{21} )-d (a_{11} )-d (a_{21} ) ]=0~\text{which implies that,}\\ \\	
				& g (x_{11}) ( [d (a_{11}+a_{21} )-d (a_{11} )-d (a_{21} ) ]_{11} \\&+ [d (a_{11}+a_{21} )-d (a_{11} ) -d (a_{21} ) ]_{21} )=0\\
					\text{and}~	& g (x_{11} ) ( [d (a_{11}+a_{21} )-d (a_{11} )-d (a_{21} ) ]_{12} \\& +[d (a_{11}+a_{21} )-d (a_{11} )-d (a_{21} ) ]_{22})=0.
			\end{aligned}
		\end{equation}
		
		Similarly, taking $x_{22} \in R_{22}$, we have
		\begin{equation}
			\label{eq2.8}
			\begin{aligned}
				& g (x_{22} ) [d (a_{11}+a_{21} )-d (a_{11} )-d (a_{21} ) ]=0 \text{ which implies that,}\\ \\
				&  g (x_{22} ) ( [d (a_{11}+a_{21} )-d (a_{11} )-d (a_{21} ) ]_{11}\\& + [d (a_{11}+a_{21} )-d (a_{11} )-d (a_{21} ) ]_{21} )=0 \\
				\text{and}~	&  g (x_{22} ) ( [d (a_{11}+a_{21} )-d (a_{11} )-d (a_{21} ) ]_{12} \\& + [d (a_{11}+a_{21} ) -d (a_{11} )-d (a_{21} ) ]_{22} )=0.
			\end{aligned}
		\end{equation}
		
		Using the Assumption \ref{asm1.1} (ii), \eqref{eq2.7} and \eqref{eq2.8},
		$$
		\begin{aligned}
			&{ [d (a_{11}+a_{21} )-d (a_{11} )-d (a_{21} ) ]_{1 k}=0} \\
			& \&~[d (a_{11}+a_{21} )-d (a_{11} )-d (a_{21} ) ]_{2 k}=0.
		\end{aligned}
		$$
		for $k=1,2$. Hence, the proof of (ii) of Lemma \ref{lem2.3}. Similarly, we can prove (iii), (iv) and (v) of Lemma \ref{lem2.3}.
	\end{proof}

	\begin{lemma}
		\label{lem2.4}
		Let $a_{ij}, b_{ij}, c_{ij} \in R_{ij}$. Then
		$$
		\begin{aligned}
			(i)~ d (a_{12}+b_{12} c_{22} )=d (a_{12} )+d (b_{12} c_{22} ),\\
			(ii)~ d (a_{21}+b_{22} c_{21} )=d (a_{21} )+d (b_{22} c_{21} ).
		\end{aligned}
		$$
	\end{lemma}

	\begin{proof} Note that, $a_{12}+b_{12} c_{22} = (e_1+b_{12} ) (a_{12}+c_{22} )$. Hence,
		$$
		\begin{aligned}
			d (a_{12}+b_{12} c_{22} ) &=d [ (e_1+b_{12} ) (a_{12}+c_{22} ) ]\\
			&=d (e_1+b_{12} ) \alpha (a_{12}+c_{22} )+g (e_1+b_{12} ) d (a_{12}+c_{22} )\\
			&=d (e_1 ) \alpha (a_{12} )+g (e_1 ) d (a_{12} )+d (b_{12} ) \alpha (a_{12})  +g (b_{12} ) d (a_{12} )\\
			&+d (e_1 ) \alpha (c_{22} )+g (e_1 ) d (c_{22} )+d (b_{12} ) \alpha (c_{22} )+g (b_{12} ) d (c_{22} )\\
			&\text{(By the Assumption \ref{asm1.1} (iii), Lemma \ref{lem2.3} (i) and (iii))}\\
			&=d (e_1 a_{12} )+d (b_{12} a_{12} )+d (e_1 c_{22} )+d (b_{12} c_{22} )	\\
			&=d (a_{12} )+d (b_{12}c_{22}),~\text{By Lemma \ref{lem2.2}}.
		\end{aligned}
		$$
		Hence, the proof of (i) of Lemma \ref{lem2.4}. Similarly, to prove (ii) of Lemma \ref{lem2.4}, we use the identity
		$a_{21}+b_{22} c_{21}= (a_{21}+b_{22} ) (e_1+c_{21} ).$ Then using the Assumption \ref{asm1.1} (iv), Lemma \ref{lem2.3} (ii) and (iv), we have the desired result.
	\end{proof}
	
	\begin{lemma}
		\label{lem2.5}
		Let $a_{ij}, b_{ij} \in R_{ij}$. Then
		$$
		\begin{aligned}
			&{ (i)~ } d (a_{12}+b_{12} )=d (a_{12} )+d (b_{12} ), \\
			& { (ii)~ } d (a_{21}+b_{21} )=d (a_{21} )+d (b_{21} ).
		\end{aligned}
		$$
	\end{lemma}
	
	\begin{proof}
		Let $x_{22} \in R_{22}$. Then
		\begin{equation}
			\label{eq2.9}
			\begin{aligned}
				&d (a_{12}+b_{12} ) g (x_{22} )+\alpha (a_{12}+b_{12} ) d (x_{22} )\\
				&=d ( (a_{12}+b_{12} ) x_{22} ) \\
				&=d ( a_{12}x_{22} + b_{12} x_{22} )\\
				&=d (a_{12} x_{22} )+d (b_{12} x_{22} )~\text{(By Lemma \ref{lem2.4} (i))} \\
				&=d (a_{12} ) g (x_{22} )+\alpha (a_{12} ) d (x_{22} )+d (b_{12} ) g (x_{22} )+\alpha (b_{12} ) d (x_{22})
			\end{aligned}
		\end{equation}
		
		Comparing both sides of \eqref{eq2.9},
		\begin{equation}
			\label{eq2.10}
			\begin{aligned}
				&{ [d (a_{12}+b_{12} )-d (a_{12} )-d (b_{12} ) ] g (x_{22} )=0  }~\text{which implies that} \\ \\
				&  ( [d (a_{12}+b_{12} )-d (a_{12} )-d (b_{12} ) ]_{11} \\&+ { [d (a_{12}+b_{12} )-d (a_{12} ) }-d (b_{12} ) ]_{12} ) g (x_{22} )=0 \\
				\text{and}~	&  ( [d (a_{12}+b_{12} )-d (a_{12} )-d (b_{12} ) ]_{21} \\&+ { [d (a_{12}+b_{12} )-d (a_{12} ) }-d (b_{12} ) ]_{22} ) g (x_{22} )=0.
			\end{aligned}
		\end{equation}
		
		Let $x_{11} \in R_{11}$. Then
		\begin{equation}
			\label{eq2.11}
			\begin{aligned}
				&d (a_{12}+b_{12} ) g (x_{11} )+\alpha (a_{12}+b_{12} ) d (x_{11} )\\
				&=d ( (a_{12}+b_{12} ) x_{11} ) \\
				&=d(0)=0~\text{(By Lemma \ref{lem2.2})} \\
				&=d (a_{12} x_{11} )+d (b_{12} x_{11} ) \\
				&=d (a_{12} ) g (x_{11} )+\alpha (a_{12} ) d (x_{11} )+d (b_{12} ) g (x_{11} )+\alpha (b_{12} ) d (x_{11} ).
			\end{aligned}
		\end{equation}
		
		Comparing both sides of \eqref{eq2.11},
		\begin{equation}
			\label{eq2.12}
			\begin{aligned}
				&{ [d (a_{12}+b_{12} )-d (a_{12} )-d(b_{12}) ] g (x_{11} )=0 } ~\text{which implies that,} \\ \\
				&  ( [d (a_{12}+b_{12} )-d (a_{12} )-d (b_{12} ) ]_{11} \\&+  { [d (a_{12}+b_{12} )-d (a_{12} ) }-d (b_{12} ) ]_{12} ) g (x_{11} )=0 \\
				\text{and}~	& ( [d (a_{12}+b_{12} )-d(a_{12})-d (b_{12} ) ]_{21} \\& + { [d (a_{12}+b_{12} )-d (a_{12} ) }-d (b_{12} ) ]_{22} ) g (x_{11} )=0.
			\end{aligned}
		\end{equation}
		
		By the Assumption \ref{asm1.1} (i), \eqref{eq2.10} and \eqref{eq2.12},
		$$
		\begin{aligned}
			& { [d (a_{12}+b_{12} )-d (a_{12} )-d (b_{12} ) ]_{k 2}=0 } \\
			\& \quad & { [d (a_{12}+b_{12} )-d (a_{12} )-d (b_{12} ) ]_{k 1}=0. }
		\end{aligned}
		$$
		for $k=1,2$. Hence, we have
		\begin{align*}
			d (a_{12}+b_{12} )=d (a_{12} )+d (b_{12} ).
		\end{align*}	
		
		Similarly, we can prove that
		\begin{align*}
			d (a_{21}+b_{21} )=d (a_{21} )+d (b_{21}).
		\end{align*}	
	\end{proof}
	
	\begin{lemma}
		\label{lem2.6}
		Let $a_{ii}, b_{ii} \in R_{ii}$. Then
		$$
		\begin{aligned}
			&{ (i)~ } d (a_{11}+b_{11} )=d (a_{11} )+d (b_{11} ), \\
			& { (ii)~ } d (a_{22}+b_{22} )=d (a_{22} )+d (b_{22} ).
		\end{aligned}
		$$
	\end{lemma}
	
	\begin{proof}
		Let $x_{12} \in R_{12}$. Then
		\begin{equation}
			\label{eq2.13}
			\begin{aligned}
				&d (a_{11}+b_{11} ) g (x_{12} )+\alpha (a_{11}+b_{11} ) d (x_{12} )\\
				&= d ( (a_{11}+b_{11} ) x_{12} ) \\
				&= d (a_{11} x_{12}+b_{11} x_{12} ) \\
				&= d (a_{11} x_{12} )+d (b_{11} x_{12} )~\text{(By Lemma \ref{lem2.5} (i))} \\
				&= d (a_{11} ) g (x_{12} )+\alpha (a_{11} ) d (x_{12} )+d (b_{11} ) g (x_{12} )+\alpha (b_{11} ) d (x_{12} ).
			\end{aligned}
		\end{equation}
		
		Comparing both sides of \eqref{eq2.13},
		\begin{equation}
			\label{eq2.14}
			\begin{aligned}
				&{ [d (a_{11}+b_{11} )-d (a_{11} )-d (b_{11} ) ] g (x_{12} )=0  }~\text{which implies that,} \\ \\
				&  ( [d (a_{11}+b_{11} )-d (a_{11} )-d (b_{11} ) ]_{11} \\&+ { [d (a_{11}+b_{11} )-d (a_{11} ) }-d (b_{11} ) ]_{12} ) g (x_{12} )=0 \\
				\text{and}~	&  ( [d (a_{11}+b_{11} )-d (a_{11} )-d (b_{11} ) ]_{21} \\&+ { [d (a_{11}+b_{11} )-d (a_{11} ) }-d (b_{11} ) ]_{22} ) g (x_{12} )=0.
			\end{aligned}
		\end{equation}

		Similarly, taking $x_{22} \in R_{22}$, we have
		\begin{equation}
			\label{eq2.15}
			\begin{aligned}
				&{ [d (a_{11}+b_{11} )-d (a_{11} )-d (b_{11} ) ] g (x_{22} )=0  }~\text{which implies that,} \\ \\
				&  ( [d (a_{11}+b_{11} )-d (a_{11} )-d (b_{11} ) ]_{11} \\&+ { [d (a_{11}+b_{11} )-d (a_{11} ) }-d (b_{11} ) ]_{12} ) g (x_{22} )=0 \\
				\text{and}~	&  ( [d (a_{11}+b_{11} )-d (a_{11} )-d (b_{11} ) ]_{21} \\&+ { [d (a_{11}+b_{11} )-d (a_{11} ) }-d (b_{11} ) ]_{22} ) g (x_{22} )=0.
			\end{aligned}
		\end{equation}
		
		By the Assumption \ref{asm1.1} (i), \eqref{eq2.14} and \eqref{eq2.15}, we have
		$$
		\begin{aligned}
			&{ [d (a_{11}+b_{11} )-d (a_{11} )-d (b_{11} ) ]_{k1}=0}, \\
			& [d (a_{11}+b_{11} )-d (a_{11} )-d (b_{11} ) ]_{k2}=0 .
		\end{aligned}
		$$
		for $k=1,2$. Hence, we have
		\begin{align*}
			d (a_{11}+b_{11} )=d (a_{11} )+d (b_{11} ).
		\end{align*}	
		
		Similarly, we can prove that
		\begin{align*}
			d (a_{22}+b_{22} )=d (a_{22} )+d (b_{22} ).
		\end{align*}	
	\end{proof}

	\begin{lemma}
		\label{lem2.7}
		Let $a_{ij}\in R_{ij}$. Then
		\begin{align*}
			d (a_{11}+a_{12}+a_{21}+a_{22} )=d (a_{11} )+d (a_{12} )+d (a_{21} )
			+d (a_{22} ).
		\end{align*}
	\end{lemma}
	
	\begin{proof}
		Let $x_{11}\in R_{11}$. Then
		\begin{equation}
			\label{eq2.16}
			\begin{aligned}
				&d (a_{11}+a_{12}+a_{21}+a_{22} ) g (x_{11} )
				+\alpha (a_{11}+a_{12}+a_{21}+a_{22} ) d (x_{11} )\\
				&= d ((a_{11}+a_{12}+a_{21}+a_{22} ) x_{11} ) \\
				&=d (a_{11} x_{11}+a_{21} x_{11} ) \\
				&=d (a_{11} x_{11} )+d (a_{12} x_{11} )+d (a_{21} x_{11} )+d (a_{22} x_{11} ) \\
				&=d (a_{11} ) g (x_{11} )+\alpha (a_{11} ) d (x_{11} )+d (a_{12} ) g (x_{11} )+\alpha (a_{12} ) d (x_{11} ) \\
				&+d (a_{21} ) g (x_{11} )+\alpha (a_{21} ) d (x_{11} )+d (a_{22} ) g (x_{11} ) +\alpha (a_{22} ) d (x_{11} ).
			\end{aligned}
		\end{equation}
		
		Comparing both sides of \eqref{eq2.16},
		\begin{equation}
			\label{eq2.17}
			\begin{aligned}
				& [d (a_{11}+a_{12}+a_{21}+a_{22} )-d (a_{11} )-d (a_{12} ) \\&-d (a_{21}) -d (a_{22} ) ] g (x_{11} )=0 \\
				& ~\text{which implies that,}\\
				&([d (a_{11}+a_{12}+a_{21}+a_{22} )-d (a_{11} )-d (a_{12} ) -d (a_{21}) -d (a_{22} ) ]_{11}\\
				&+[d (a_{11}+a_{12}+a_{21}+a_{22} )-d (a_{11} )-d (a_{12} )\\&-d (a_{21}) -d (a_{22} ) ]_{12})g (x_{11} )=0\\
				\text{and}~	&([d (a_{11}+a_{12}+a_{21}+a_{22} )-d (a_{11} )-d (a_{12} )-d (a_{21}) -d (a_{22} ) ]_{21}\\
				&+[d (a_{11}+a_{12}+a_{21}+a_{22} )-d (a_{11} )\\&-d (a_{12} )-d (a_{21}) -d (a_{22} ) ]_{22})g (x_{11} )=0.
			\end{aligned}
		\end{equation}
		
		Similarly, taking $x_{22} \in R_{22}$, we have
		\begin{equation}
			\label{eq2.18}
			\begin{aligned}
				& [d (a_{11}+a_{12}+a_{21}+a_{22} )-d (a_{11} )-d (a_{12} ) \\&-d (a_{21}) -d (a_{22} ) ] g (x_{22} )=0  ~\text{which implies that,}\\ \\
				&([d (a_{11}+a_{12}+a_{21}+a_{22} )-d (a_{11} )-d (a_{12} )-d (a_{21}) -d (a_{22} ) ]_{11}\\
				&+[d (a_{11}+a_{12}+a_{21}+a_{22} )-d (a_{11} )-d (a_{12} ) \\&-d (a_{21}) -d (a_{22} ) ]_{12})g (x_{22} )=0\\
				\text{and}~	&([d (a_{11}+a_{12}+a_{21}+a_{22} )-d (a_{11} )-d (a_{12} )-d (a_{21}) -d (a_{22} ) ]_{21}\\
				&+[d (a_{11}+a_{12}+a_{21}+a_{22} )-d (a_{11} )-d (a_{12} ) \\&-d (a_{21}) -d (a_{22} ) ]_{22})g (x_{22} )=0.
			\end{aligned}
		\end{equation}
		
		By the Assumption \ref{asm1.1} (i), \eqref{eq2.17} and \eqref{eq2.18}, we get
		\begin{align*}
			& [d (a_{11}+a_{12}+a_{21}+a_{22} )-d (a_{11} )-d (a_{12} )-d (a_{21}) -d (a_{22} ) ]_{k1}=0,\\
			& [d (a_{11}+a_{12}+a_{21}+a_{22} )-d (a_{11} )-d (a_{12} )-d (a_{21}) -d (a_{22} ) ]_{k2}=0,
		\end{align*}
		for $k=1,2$. Hence,
		\begin{align*}
			d (a_{11}+a_{12}+a_{21}+a_{22} )=d (a_{11} )+d (a_{12} )+d (a_{21} )
			+d (a_{22} ).
		\end{align*}
	\end{proof}

	\begin{proof} [Proof of Theorem \ref{thm2.1}] Let $a,b \in R$. Then
		\begin{align*}
			&a=a_{11}+a_{12}+a_{21}+a_{22},\\
			& b= b_{11}+b_{12}+b_{21}+b_{22},
		\end{align*}
		for some $a_{ij}, b_{ij} \in R_{ij}$. Now,
		$$
		\text {  } \begin{aligned}
			d(a+b)=~& d (a_{11}+a_{12}+a_{21}+a_{22}+b_{11}+b_{12}+b_{21}+b_{22} ) \\
			=~& d ( (a_{11}+b_{11} )+ (a_{12}+b_{12} )+ (a_{21}+b_{21} )+ (a_{22}+b_{22} ) ) \\
			=~& d (a_{11}+b_{11} )+d (a_{12}+b_{12} )+d (a_{21}+b_{21} ) \\&+d (a_{22}+b_{22} ) ~\text{(By Lemma \ref{lem2.7})}\\
			=~& d (a_{11} )+d (b_{11} )+d (a_{12} )+d (b_{12} )+d (a_{21} )\\& +d (b_{21} )+d (a_{22} )+d (b_{22} ) ~\text{(By Lemma \ref{lem2.5} and \ref{lem2.6})}\\
			=~& d (a_{11}+a_{12}+a_{21}+a_{22} )+d (b_{11}+b_{12}+b_{21}+b_{22} ) \\
			&~\text{(By Lemma \ref{lem2.7})}\\
			=~& d(a)+d(b).
		\end{aligned}
		$$
		Hence, $d$ is additive on $R$. Moreover, $d$ is a skew semi-derivation on $R$.
	\end{proof}
	
\begin{cor}
\label{cor2.1}
Let $R$ be a ring containing a non-trivial idempotent $e$. Let $d$ be a multiplicative skew semi-derivation on $R$, and $g$ be a map on $R$ such that $g(0)=0$. Suppose $d$ and $g$ satisfy the following conditions for $i, j, k=1,2$:\\
	(i) If $(a_{k 1}+b_{k 2}) g\ (x_{i j})=0$ and $i \leq j$ for all $x_{i j} \in R_{i j}$, then
	$$
	\begin{cases}a_{k 1}=0 & \text { if } i=1, \\ b_{k 2}=0 & \text { if } i=2 .\end{cases}
	$$
	(ii) If $g (x_{i i})\ (a_{1 j}+b_{2 j})=0$ for all $x_{i i} \in R_{i i}$, then
	$$
	\begin{cases}a_{1 j}=0 & \text { if } i=1, \\ b_{2 j}=0 & \text { if } i=2 .\end{cases}
	$$
	(iii) If $d (x_{i i}) (a_{1 j}+b_{2 j})=0$ for all $x_{i i} \in R_{i i}$, then
	$$
	\begin{cases}a_{1 j}=0 & \text { if } i=1, \\ b_{2 j}=0 & \text { if } i=2 .\end{cases}
	$$
	Then $d$ is additive.
\end{cor}

\begin{proof}
It is sufficient to prove that if condition (iii) of the above corollary holds, then it implies conditions (iii) and (iv) of Assumption \ref{asm1.1}. We prove conditions (iii) and (iv) of Assumption \ref{asm1.1}.
	Let $x_{11}, a_{11} \in R_{11}$ and $b_{12} \in R_{12}$.
	 Then
	$$
	d (x_{11}) g (a_{11}+b_{12})=d (x_{11} (a_{11}+b_{12}))- \alpha(x_{11}) d (a_{11}+b_{12}).
	$$
By Lemma \ref{lem2.3}(i),
	$$
	\begin{aligned}
		d (x_{11}) g (a_{11}+b_{12}) & =d (x_{11} a_{11})+d (x_{11} b_{12})- \alpha(x_{11}) d (a_{11})- \alpha(x_{11}) d (b_{12}) \\
		& =d (x_{11}) g (a_{11})+d (x_{11}) g (b_{12}),
	\end{aligned}
	$$
	implies that
	$$
	d (x_{11}) [g (a_{11}+b_{12})-g (a_{11})-g (b_{12})]=0.
	$$
Consequently,
	$$
	d (x_{11}) ( [g (a_{11}+b_{12})-g (a_{11})-g (b_{12})]_{11}+ [g (a_{11}+b_{12})-g (a_{11})-g (b_{12})]_{21})=0
	$$
	and
	$$
	d (x_{11}) ( [g (a_{11}+b_{12})-g (a_{11})-g (b_{12})]_{12}+ [g (a_{11}+b_{12})-g (a_{11})-g (b_{12})]_{22})=0.
	$$
Similarly, taking $x_{22} \in R_{22}$, we have
	$$
	d (x_{22}) ( [g (a_{11}+b_{12})-g (a_{11})-g (b_{12})]_{11}+ [g (a_{11}+b_{12})-g (a_{11})-g (b_{12})]_{21})=0
	$$
	and
	$$
	d(x_{22}) ( [g (a_{11}+b_{12})-g(a_{11})-g (b_{12})]_{12}+ [g (a_{11}+b_{12})-g (a_{11})-g (b_{12})]_{22})=0 .
	$$
	Using condition (iii) of corollary \ref{cor2.1},
	$$
	\begin{aligned}
		& { [g(a_{11}+b_{12})-g (a_{11})-g (b_{12})]_{1 k}=0,} \\
		& { [g (a_{11}+b_{12})-g (a_{11})-g (b_{12})]_{2 k}=0,}
	\end{aligned}
	$$
	where $k=1,2$. Hence,
	\begin{align*}
	g (a_{11}+b_{12})=g (a_{11})+g (b_{12}).
\end{align*}	
Similarly, we can prove that
\begin{align*}
	g (a_{11}+b_{21})=g(a_{11})+g (b_{21}).
\end{align*}	
\end{proof}
\begin{cor}
\label{cor2.2}
Let $R$ be a prime ring with identity $1\neq 0$ and containing an idempotent $e \neq 0,1$. If $d$ is any multiplicative skew semi-derivation on $R$ where the associated map $g$ is an endomorphism on $R$ such that $g(e)=e$, then $d$ is additive.
\end{cor}

	\begin{cor}
	Let $R$ be a ring as in Corollary \ref{cor2.2}. If $d$ is any multiplicative skew semi-derivation on $R$ where the associated map $g$ is invariant on $R_{i j}$, i.e., $g (R_{i j}) \subseteq R_{i j}$ where $i, j=1,2$ and $g$ satisfies conditions (iii) and (iv) of Assumption \ref{asm1.1}, then $d$ is additive.
		\end{cor}
	
\begin{theorem}
\label{thm2.11}
Let $R$ be a ring with a non-trivial idempotent $e$ and $g$ is any map on $R$ which vanishes at zero and satisfies the following conditions, for $i, j=1,2$:
\begin{align*}
& (i) ~\text{If}~  (a_{j1}+b_{j2}) g (x_{ii})=0 ~\text{for all } x_{ii} \in R_{ii}, \text{ then } a_{j 1}+b_{j 2}=0.\\
	&(ii)~ \text{If}~ g (x_{i i})(a_{11}+a_{21})=0 ~\text{ for all } x_{i i} \in R_{i i}, \text{ then } a_{11}+a_{21}=0.\\
\end{align*}
	If $d$ is a multiplicative skew semi-derivation on $R$ with the associated map $g$ and automorphism $\alpha$, then $d$ is additive.
\end{theorem}
	
	To prove Theorem \ref{thm2.11}, we prove some results as follows.
	Note that, we have $d(0)=0$ by Lemma \ref{lem2.2}.
	
\begin{lemma}
\label{lem2.12}
Let $a_{11} \in R_{11},~ b_{12} \in R_{12}$ and $c_{21} \in R_{21}$. Then
\begin{align*}
&(i)~ d (a_{11}+b_{12})=d (a_{11})+d (b_{12}),\\
&	(ii)~ d (a_{11}+c_{21})=d(a_{11})+d (c_{21}).
\end{align*}
\end{lemma}

\begin{proof}
 (i) Let $a_{11}, x_{11} \in R_{11}$ and $b_{12} \in R_{12}$. Then
	$$
	\begin{aligned}
		d ( (a_{11}+b_{12}) x_{11}) & =d (a_{11} x_{11}) \\
		& =d (a_{11} x_{11})+d (b_{12} x_{11}) \\
		& =d (a_{11}) g (x_{11})+ \alpha(a_{11}) d (x_{11}) \\&+d (b_{12}) g (x_{11})+ \alpha(b_{12}) d (x_{11}) .
	\end{aligned}
	$$
	On the other hand,
	$$
	d ( (a_{11}+b_{12}) x_{11})=d (a_{11}+b_{12}) g (x_{11})+ \alpha(a_{11}+b_{12}) d (x_{11}) .
	$$
	Comparing the identities,
	$$
	[d (a_{11}+b_{12})-d (a_{11})-d (b_{12})] g (x_{11})=0
	$$
	which implies that
	$$
	d\cdot g (x_{11})=0,
	$$
	where $d=d (a_{11}+b_{12})-d (a_{11})-d (b_{12})$. This gives
	\begin{equation}
	\label{eq2.19a}
	([d]_{11}+[d]_{12}) g (x_{11})=0
	\end{equation}
	and
	\begin{equation}
	\label{eq2.20a}
	([d]_{21}+[d]_{22}) g (x_{11})=0.
	\end{equation}
	Similarly, taking $x_{22} \in R_{22}$, we have
	$$
	\begin{aligned}
		d (x_{22} (a_{11}+b_{12})) & =d (x_{22} a_{11})+d (x_{22} b_{12}) \\
		& =d (x_{22}) \alpha(a_{11})+g (x_{22}) d (a_{11}) \\&+d (x_{22}) \alpha(b_{12})+g (x_{22}) d (b_{12}).
	\end{aligned}
	$$
	On the other hand,
	$$
	d (x_{22} (a_{11}+b_{12}))=d (x_{22}) \alpha(a_{11}+b_{12})+g (x_{22}) d (a_{11}+b_{12}).
	$$
	Comparing,
	$$
	g (x_{22}) [d (a_{11}+b_{12})-d (a_{11})-d (b_{12})]=0
	$$
	which implies
	$$
	g (x_{22})\cdot d=0.
	$$
	We get
	\begin{equation}
	\label{eq2.21a}
	g (x_{22}) ([d]_{11}+[d]_{21})=0 .
	\end{equation}
	Using Theorem \ref{thm2.11} $(i)-(ii)$ in \eqref{eq2.19a}-\eqref{eq2.21a}, we have
	$$
	\begin{aligned}
		& {[d]_{11}+[d]_{12}=0,} \\
		& {[d]_{21}+[d]_{22}=0,} \\
		& {[d]_{11}+[d]_{21}=0.}
	\end{aligned}
	$$
	Therefore, $[d]_{i j}=0$ for all $i$, $j$, which completes the proof.\\
	
	(ii) Let $a_{11} \in R_{11}, x_{22} \in R_{22}$ and $c_{21} \in R_{21}$. Then
	\begin{align*}
	d ( (a_{11}+c_{21}) x_{22})&=d (a_{11} x_{22})+d (c_{21} x_{22})\\
	&=d (a_{11}) g (x_{22})+ \alpha(a_{11}) d (x_{22}) \\&+d (c_{21}) g (x_{22})+ \alpha(c_{21}) d (x_{22}).
	\end{align*}
	
	On the other hand,
	\begin{align*}
	d( (a_{11}+c_{21}) x_{22})=d (a_{11}+c_{21}) g (x_{22})+ \alpha(a_{11}+c_{21}) d (x_{22}).
	\end{align*}
	
	Comparing these two equalities, we have
	$$
	 [d (a_{11}+c_{21})-d (a_{11})-d (c_{21})] g (x_{22})=0.
	$$
	This implies
	$$
	d \cdot g (x_{22})=0
	$$
	where $d=d (a_{11}+c_{21})-d (a_{11})-d (c_{21})$. Now, we get
	\begin{equation}
	\label{eq2.19}
	([d]_{11}+[d]_{12}) g (x_{22})=0
	\end{equation}
	and
\begin{equation}
\label{eq2.20}
	([d]_{21}+[d]_{22}) g (x_{22})=0.
\end{equation}
	
	Similarly, taking $x_{11} \in R_{11}$, we have
	$$
	\begin{aligned}
		d (x_{11} (a_{11}+c_{21})) & =d (x_{11} a_{11}) \\
		& =d (x_{11} a_{11})+d (x_{11} c_{21}) \\
		& =d (x_{11}) \alpha(a_{11})+g (x_{11}) d (a_{11}) \\&+d (x_{11}) \alpha(c_{21})+g (x_{11}) d (c_{21}).
	\end{aligned}
	$$
	On the other hand,
	$$
	d (x_{11} (a_{11}+c_{21}))=d (x_{11}) \alpha(a_{11}+c_{21})+g (x_{11}) d (a_{11}+c_{21})
	$$
	Comparing the above equalities, we get
	$$
	g (x_{11}) [d (a_{11}+c_{21})-d (a_{11})-d (c_{21})]=0
	$$
	which implies that
	$$
	g (x_{11})\cdot d=0.
	$$
	Also, we have
	\begin{equation}
	\label{eq2.21}
	g (x_{11}) ([d]_{11}+[d]_{21})=0.
	\end{equation}
	
	Using Theorem \ref{thm2.11} (i)-(ii) in \eqref{eq2.19}-\eqref{eq2.21},
	$$
	\begin{aligned}
		& {[d]_{11}+[d]_{12}=0} \\
		& {[d]_{21}+[d]_{22}=0} \\
		& {[d]_{11}+[d]_{21}=0}
	\end{aligned}
	$$
	Therefore, $[d]_{i j}=0$ for all $i, j$.
	Hence, we have the desired result.
\end{proof}

\begin{lemma}
	\label{lem2.13}
The following statements hold:\\
	(i) $d$ is additive on $R_{12}$,\\
	(ii) $d$ is additive on $R_{21}$,\\
	(iii) $d$ is additive on $R_{11}$,\\
	(iv) $d$ is additive on $R_{22}$.
\end{lemma}

\begin{proof}
Let $x_{11} \in R_{11}$ and $a_{12},~ b_{12} \in R_{12}$. Then
	$$
	\begin{aligned}
		d ( (a_{12}+b_{12}) x_{11}) & =d (a_{12} x_{11})+d (b_{12} x_{11}) \\
		& =d (a_{12}) g (x_{11})+ \alpha(a_{12}) d (x_{11}) \\&+d (b_{12}) g (x_{11})+ \alpha(b_{12}) d (_{11}).
	\end{aligned}
	$$
On the other hand,
$$
d( (a_{12}+b_{12}) x_{11})=d (a_{12}+b_{12}) g (x_{11})+ \alpha(a_{12}+b_{12}) d (x_{11}) .
$$
Comparing above identities,
$$
 [d (a_{12}+b_{12})-d (a_{12})-d (b_{12})] g (x_{11})=0
$$
which implies,
$$
d \cdot g (x_{11})=0,
$$
where $d=d (a_{12}+b_{12})-d (a_{12})-d (b_{12})$. Now, we get
\begin{equation}
\label{eq2.22}
([d]_{11}+[d]_{12}) g (x_{11})=0
\end{equation}
and
\begin{equation}
\label{eq2.23}
([d]_{21}+[d]_{22}) g (x_{11})=0.
\end{equation}
Similarly, taking $x_{22} \in R_{22}$,
$$
\begin{aligned}
	d (x_{22} (a_{12}+b_{12})) & =d (x_{22} a_{12})+d(x_{22} b_{12}) \\
	& =d (x_{22}) \alpha(a_{12})+g (x_{22}) d (a_{12}) \\&+d (x_{22}) \alpha(b_{12})+g (x_{22}) d (b_{12}) .
\end{aligned}
$$
On the other hand,
$$
d(x_{22} (a_{12}+b_{12}))=d (x_{22}) \alpha(a_{12}+b_{12})+g (x_{22}) d (a_{12}+b_{12}) .
$$
Comparing the above equalities,
$$
g (x_{22}) [d (a_{12}+b_{12})-d (a_{12})-d (b_{12})]=0
$$
which implies,
$$
g(x_{22})\cdot d=0 .
$$
We get
\begin{equation}
\label{eq2.24}
g (x_{22}) ([d]_{11}+[d]_{21})=0.
\end{equation}
Using Theorem \ref{thm2.11} (i)-(ii) in \eqref{eq2.22}-\eqref{eq2.24},
$$
\begin{aligned}
	& {[d]_{11}+[d]_{12}=0,} \\
	& {[d]_{21}+[d]_{22}=0,} \\
	& {[d]_{11}+[d]_{21}=0 .}
\end{aligned}
$$
Therefore, we get $[d]_{i j}=0$ for all $i, j$. This completes the proof of $(i)$ of Lemma \ref{lem2.13}. Similarly, we can prove $(ii), (iii)$ and $(iv)$ of Lemma \ref{lem2.13}.
\end{proof}

\begin{lemma}
	\label{lem2.14}
 Let $a_{ij} \in R_{ij}$. Then
  $$ d (a_{11}+a_{12}+a_{21}+a_{22})=d (a_{11})+d (a_{12})+d (a_{21})+d (a_{22}).$$
\end{lemma}

\begin{proof}
 Let $x_{11} \in R_{11}$. Then
\begin{align*}
&d ( (a_{11}+a_{12}+a_{21}+a_{22}) x_{11})\\
&=d (a_{11}+a_{12}+a_{21}+a_{22}) g (x_{11})+\alpha(a_{11}+a_{12}+a_{21}+a_{22}) d (x_{11}).
\end{align*}
On the other hand, by Lemma \ref{lem2.12} $(ii)$, we have
$$
\begin{aligned}
	&d ( (a_{11}+a_{12}+a_{21}+a_{22}) x_{11})\\
	=&  d (a_{11} x_{11}+a_{21} x_{11}) \\
	= & d (a_{11} x_{11})+d (a_{21} x_{11}) \\
	= & d (a_{11} x_{11})+d(a_{12} x_{11})+d (a_{21} x_{11})+d (a_{22} x_{11}) \\
	= & d (a_{11}) g(x_{11})+ \alpha(a_{11}) d (x_{11})+d (a_{12}) g (x_{11})+ \alpha(a_{12}) d (x_{11}) \\
	& +d (a_{21}) g (x_{11})+ \alpha(a_{21}) d (x_{11})+d (a_{22}) g (x_{11})+ \alpha(a_{22}) d (x_{11}).
\end{aligned}
$$
Comparing the above equalities, we get
$$
[d (a_{11}+a_{12}+a_{21}+a_{22})-d(a_{11})-d (a_{12})-d(a_{21})-d (a_{22})] g (x_{11})=0,
$$
which implies that
$$
d \cdot g (x_{11})=0,
$$
where $d=d (a_{11}+a_{12}+a_{21}+a_{22})-d (a_{11})-d (a_{12})-d (a_{21})-d (a_{22})$. Therefore,
\begin{equation}
\label{eq2.25}
([d]_{11}+[d]_{12}) g (x_{11})=0
\end{equation}
and
\begin{equation}
\label{eq2.26}
 ([d]_{21}+[d]_{22}) g (x_{11})=0.
\end{equation}
Similarly,
\begin{align*}
&d (x_{11} (a_{11}+a_{12}+a_{21}+a_{22}))\\
&=d (x_{11})\alpha (a_{11}+a_{12}+a_{21}+a_{22})+g (x_{11})d(a_{11}+a_{12}+a_{21}+a_{22}).
\end{align*}
On the other hand, by Lemma \ref{lem2.12} $(i)$, we have
$$
\begin{aligned}
	&d (x_{11} (a_{11}+a_{12}+a_{21}+a_{22}))\\
	= & d (x_{11} a_{11}+x_{11} a_{12}) \\
	= & d (x_{11} a_{11})+d (x_{11} a_{12}) \\
	= & d (x_{11} a_{11})+d (x_{11} a_{12})+d(x_{11} a_{21})+d (x_{11} a_{22}) \\
	= & d (x_{11}) \alpha(a_{11})+g(x_{11}) d (a_{11})+d (x_{11}) \alpha(a_{12})+g (x_{11}) d (a_{12}) \\
	& +d(x_{11}) \alpha(a_{21})+g (x_{11}) d (a_{21})+d (x_{11}) \alpha(a_{22})+g (x_{11})  d(a_{22}).
\end{aligned}
$$
Comparing the above identities, we get
$$
g (x_{11}) [d (a_{11}+a_{12}+a_{21}+a_{22})-d (a_{11})-d (a_{12})-d (a_{21})-d (a_{22})]=0
$$
which implies
$$
g (x_{11})\cdot d=0,
$$
where $d=d (a_{11}+a_{12}+a_{21}+a_{22})-d (a_{11})-d (a_{12})-d (a_{21})-d (a_{22})$. Also,
\begin{equation}
\label{eq2.27}
g (x_{11}) ([d]_{11}+[d]_{21})=0.
\end{equation}
Using Theorem \ref{thm2.11} $(i)-(ii)$ in \eqref{eq2.25}-\eqref{eq2.27}
$$
\begin{aligned}
	& {[d]_{11}+[d]_{12}=0} \\
	& {[d]_{21}+[d]_{22}=0} \\
	& {[d]_{11}+[d]_{21}=0}.
\end{aligned}
$$
Therefore, we have $[d]_{i j}=0$ for all $i, j$. 
\end{proof}

\begin{proof} [Proof of Theorem \ref{thm2.11}]
	Let $a, b \in R$. Then $a=a_{11}+a_{12}+a_{21}+a_{22}$ and $b=b_{11}+b_{12}+b_{21}+b_{22}$ for $a_{ij},~b_{ij} \in R_{ij}$.
	Using Lemma \ref{lem2.13} and \ref{lem2.14},
$$
\begin{aligned}
	d(a+b) & =d (a_{11}+a_{12}+a_{21}+a_{22}+b_{11}+b_{12}+b_{21}+b_{22}) \\
	& =d ( (a_{11}+b_{11})+ (a_{12}+b_{12})+ (a_{21}+b_{21})+ (a_{22}+b_{22}))\\
	&=d (a_{11}+b_{11})+d (a_{12}+b_{12})+d (a_{21}+b_{21})+d(a_{22}+b_{22})\\
    &=d (a_{11})+d (b_{11})+d (a_{12})+d (b_{12})+d (a_{21}) \\& ~+d (b_{21})+d (a_{22})+d (b_{22})\\
   & =d(a_{11}+a_{12}+a_{21}+a_{22})+d (b_{11}+b_{12}+b_{21}+b_{22}) \\
	& =d(a)+d(b).
\end{aligned}
$$
Hence $d$ is additive.
\end{proof}

\begin{cor}
Let $R$ be a ring containing a nontrivial idempotent $e$ and $d$ be a multiplicative skew semi-derivation on $R$. Let $g$ be a map on $R$ with $g(0)=0$ and satisfies three of the following conditions for $i=1,2~:$ \\

(i) If $ (a_{11}+b_{12}) g (x_{ii})=0$ for all $x_{ii} \in R_{ii}$, then $a_{11}+b_{12}=0$.\\

(ii) If $ (a_{21}+b_{22}) g (x_{ii})=0$ for all $x_{ii} \in R_{ii}$, then $a_{21}+b_{22}=0$.\\

(iii) If $g (x_{ii}) (a_{11}+b_{21})=0$ for all $x_{ii} \in R_{ii}$, then $a_{11}+b_{21}=0$.\\

(iv) If $g(x_{ii}) (a_{12}+b_{22})=0$ for all $x_{ii} \in R_{ii}$, then $a_{12}+b_{22}=0$.\\

Then $d$ is additive.	
\end{cor}

	\section{Generalized Skew Semi-derivation}
	\begin{theorem}
		\label{thm3.1}
		Let $R$ be a ring and $ g: R \rightarrow R$ be a map that satisfies Assumption \ref{asm1.1}.
	If $ f: R \rightarrow R$ is a multiplicative generalized skew semi-derivation with the associated map $g: R \rightarrow R$, a skew semi-derivation $d$ and an automorphism $\alpha: R \rightarrow R$, then $f$ is additive on $R$.
	\end{theorem}
	
	Before proving Theorem \ref{thm3.1}, we prove several lemmas.
	
	\begin{lemma}
		\label{lem3.2}
		$f(0)=0$.
	\end{lemma}
	
	\begin{proof}
		$f(0)=f(0\cdot 0)=f(0)g(0)+\alpha(0)d(0)=0$ (By the Assumption \ref{asm1.1}, $g(0)=0$ and as $\alpha$ is an automorphism, $\alpha(0)=0$).
	\end{proof}	
	
	\begin{lemma}
		\label{lem3.3}
		Let $a_{ij} \in R_{ij}$. Then
		$$
		\begin{aligned}
			(i)~ f (a_{11}+a_{12} )=f (a_{11} )+ f(a_{12} ),\\
			(ii)~ f(a_{11}+a_{21} )= f(a_{11} )+f (a_{21} ),\\
			(iii)~ f(a_{22}+a_{12} )=f (a_{22} )+f (a_{12} ),\\
			(iv)~ f (a_{22}+a_{21} )=f (a_{22} )+f (a_{21}), \\
			(v)~ f(a_{11}+a_{22} )=f (a_{11} )+f(a_{22} ).
		\end{aligned}
		$$
	\end{lemma}
	
	\begin{proof}
		Let $x_{11} \in R_{11}$. Then
		\begin{equation}
			\label{eq3.1}
			f ( (a_{11}+a_{12} ) x_{11} )=f (a_{11}+a_{12} ) g (x_{11} )+\alpha (a_{11}+a_{12} ) d (x_{11} ).
		\end{equation}
		
		On the other hand, we have
		\begin{equation}
			\label{eq3.2}
			\begin{aligned}
				f ( (a_{11}+a_{12} ) x_{11} ) &= f (a_{11} x_{11} )+f (a_{12} x_{11} ) \\
				&=f (a_{11} ) g (x_{11} )+\alpha (a_{11} ) d (x_{11} ) \\&~+f (a_{12} ) g (x_{11} ) +\alpha (a_{12} ) d (x_{11})
			\end{aligned}
		\end{equation}
		
		Comparing \eqref{eq3.1} and \eqref{eq3.2},
		\begin{equation}
			\label{eq3.3}
			\begin{aligned}
				& { [f (a_{11}+a_{12} )-f (a_{11} )-f (a_{12} ) ] } g (x_{11} )=0 ~\text{which implies that,}\\ \\
				& ( [f (a_{11}+a_{12} )-f (a_{11} )-f (a_{12} ) ]_{11} \\&~+  { [f (a_{11}+a_{12} )-f (a_{11} ) } -f (a_{12} ) ]_{12} ) g (x_{11} )=0\\
				\text{and}~	& ( [f (a_{11}+a_{12} )-f (a_{11} )-f (a_{12} ) ]_{21} \\&~+  { [f (a_{11}+a_{12} )-f (a_{11} )} -f (a_{12} ) ]_{22} ) g (x_{11} )=0.
			\end{aligned}
		\end{equation}
		
		Similarly, for $x_{22} \in R_{22}$, we have
		\begin{equation}
			\label{eq3.4}
			\begin{aligned}
				&f (a_{11}+a_{12}) g (x_{22} )+\alpha (a_{11}+a_{12} ) d (x_{22}) \\
				&=f ( (a_{11}+a_{12} ) x_{22} )\\
				&=f (a_{12} x_{22} ) \\
				&=f (a_{11} x_{22} )+f (a_{12} x_{22} ) \\
				&=f (a_{11} ) g (x_{22} )+\alpha (a_{11} ) d (x_{22} )+f (a_{12} ) g (x_{22} )+\alpha (a_{12} ) d (x_{22} ).
			\end{aligned}
		\end{equation}	
		
		Comparing both sides of \eqref{eq3.4}, we get
		\begin{equation}
			\label{eq3.5}
			\begin{aligned}
				&\ [f (a_{11}+a_{12} )-f (a_{11} )-f (a_{12} ) ] g (x_{22} )=0~\text{which implies that,}\\ \\
				& ( [f (a_{11}+a_{12} )-f (a_{11} )-f (a_{12} ) ]_{11}\\& +[f (a_{11}+a_{12} )-f (a_{11} )-f (a_{12} ) ]_{12} ) g (x_{22} )=0\\
				\text{and}~	& ( [f (a_{11}+a_{12} )-f (a_{11} )-f (a_{12} ) ]_{21}  \\&~+  { [f (a_{11}+a_{12} )-f (a_{11} ) }-f (a_{12} ) ]_{22} ) g(x_{22} )=0 .
			\end{aligned}
		\end{equation}
		
		By the Assumption \ref{asm1.1} (i), \eqref{eq3.3} and \eqref{eq3.5},
		$$
		\begin{aligned}
			&{ [f (a_{11}+a_{12} )-f (a_{11} )-f (a_{12} ) ]_{k 1} } =0 \\
			\&~ & { [f (a_{11}+a_{12} )-f (a_{11} )-f (a_{12} ) ]_{k 2 }}=0
		\end{aligned}
		$$
		for $k=1,2$. This completes the proof of (i) of Lemma \ref{lem3.3}. \\
		
		To proof of (ii) of Lemma \ref{lem3.3}, take $x_{11} \in R_{11}$. Then
		\begin{equation}
			\label{eq3.6}
			\begin{aligned}
				& d (x_{11} ) \alpha (a_{11}+a_{21} )+g (x_{11} ) f (a_{11}+a_{21} )=f (x_{11} (a_{11}+a_{21} ) ) \\
				&=f (x_{11} a_{11} ) \\
				&=f (x_{11} a_{11} )+f (x_{11} a_{21} ) \\
				&=d (x_{11} ) \alpha (a_{11} )+g (x_{11} ) f (a_{11} )+d (x_{11} ) \alpha (a_{21} ) +g (x_{11} ) f (a_{21} ).
			\end{aligned}
		\end{equation}	
		
		Comparing both sides of \eqref{eq3.6},
		\begin{equation}
			\label{eq3.7}
			\begin{aligned}
				&g (x_{11} ) [f (a_{11}+a_{21} )-f (a_{11} )-f (a_{21} ) ]=0~\text{which implies that,}\\ \\	
				& g (x_{11}) ( [f (a_{11}+a_{21} )-f (a_{11} )-f (a_{21} ) ]_{11} \\&+ [f (a_{11}+a_{21} )-f (a_{11} ) -f (a_{21} ) ]_{21} )=0\\
				\text{and}~	& g (x_{11} ) ( [f (a_{11}+a_{21} )-f (a_{11} )-f (a_{21} ) ]_{12} \\ &~+[f (a_{11}+a_{21} )-f (a_{11} )-f (a_{21} ) ]_{22})=0.
			\end{aligned}
		\end{equation}
		
		Similarly, taking $x_{22} \in R_{22}$, we have
		\begin{equation}
			\label{eq3.8}
			\begin{aligned}
				& g (x_{22} ) [f (a_{11}+a_{21} )-f (a_{11} )-f (a_{21} ) ]=0 \text{ which implies that,}\\ \\
				&  g (x_{22} ) ( [f (a_{11}+a_{21} )-f (a_{11} )-f (a_{21} ) ]_{11}  \\&+ [f (a_{11}+a_{21} )-f (a_{11} )-f (a_{21} ) ]_{21} )=0 \\
				\text{and}~	&  g (x_{22} ) ( [f (a_{11}+a_{21} )-f (a_{11} )-f (a_{21} ) ]_{12} \\&+ [f (a_{11}+a_{21} ) -f (a_{11} )-f (a_{21} ) ]_{22} )=0.
			\end{aligned}
		\end{equation}
		
		Using the Assumption \ref{asm1.1} (ii), \eqref{eq3.7} and \eqref{eq3.8},
		$$
		\begin{aligned}
			&{ [f (a_{11}+a_{21} )-f (a_{11} )-f (a_{21} ) ]_{1 k}=0} \\
			& \&~[f (a_{11}+a_{21} )-f (a_{11} )-f (a_{21} ) ]_{2 k}=0.
		\end{aligned}
		$$
		for $k=1,2$. Hence, the proof of (ii) of Lemma \ref{lem3.3}. Similarly, we can prove (iii), (iv) and (v) of Lemma \ref{lem3.3}.
	\end{proof}

	\begin{lemma}
		\label{lem3.4}
		Let $a_{ij}, b_{ij}, c_{ij} \in R_{ij}$. Then
		$$
		\begin{aligned}
			(i)~ f (a_{12}+b_{12} c_{22} )=f (a_{12} )+f (b_{12} c_{22} ),\\
			(ii)~ f (a_{21}+b_{22} c_{21} )=f (a_{21} )+f (b_{22} c_{21} ).
		\end{aligned}
		$$
	\end{lemma}

	\begin{proof}
		Note that, $a_{12}+b_{12} c_{22} = (e_1+b_{12} ) (a_{12}+c_{22} )$. Hence,
		$$
		\begin{aligned}
			f (a_{12}+b_{12} c_{22} ) &=f [ (e_1+b_{12} ) (a_{12}+c_{22} ) ]\\
			&=f (e_1+b_{12} ) g (a_{12}+c_{22} )+\alpha (e_1+b_{12} ) d (a_{12}+c_{22} )\\
			&=f (e_1 ) g (a_{12} )+\alpha (e_1 ) d (a_{12} )+f (b_{12} ) g (a_{12})  +\alpha (b_{12} ) d (a_{12} )\\
			&+f (e_1 ) g (c_{22} )+\alpha (e_1 ) d (c_{22} )+f (b_{12} ) g (c_{22} )+\alpha (b_{12} ) d (c_{22} )\\
			&\text{(By the Assumption \ref{asm1.1} (iii), Lemma \ref{lem3.3} (i) and (iii))}\\
			&=f (e_1 a_{12} )+f (b_{12} a_{12} )+f (e_1 c_{22} )+f (b_{12} c_{22} )	\\
			&=f (a_{12} )+f (b_{12}c_{22}),~\text{By Lemma \ref{lem3.2}}.
		\end{aligned}
		$$
		Hence, the proof of (i) of Lemma \ref{lem3.4}. Similarly, to prove (ii) of Lemma \ref{lem3.4}, we use the identity
		$a_{21}+b_{22} c_{21}= (a_{21}+b_{22} ) (e_1+c_{21} ).$ Then using the Assumption \ref{asm1.1} (iv), Lemma \ref{lem3.3} (ii) and (iv), we have the desired result.
	\end{proof}
	
	\begin{lemma}
		\label{lem3.5}
		Let $a_{ij}, b_{ij} \in R_{ij}$. Then
		$$
		\begin{aligned}
			&{ (i)~ } f (a_{12}+b_{12} )=f (a_{12} )+f (b_{12} ), \\
			& { (ii)~ } f (a_{21}+b_{21} )=f (a_{21} )+f (b_{21} ).
		\end{aligned}
		$$
	\end{lemma}
	
	\begin{proof}
		Let $x_{22} \in R_{22}$. Then
		\begin{equation}
			\label{eq3.9}
			\begin{aligned}
				&f (a_{12}+b_{12} ) g (x_{22} )+\alpha (a_{12}+b_{12} ) d (x_{22} )\\
				&=f ( (a_{12}+b_{12} ) x_{22} ) \\
				&=f ( a_{12}x_{22} + b_{12} x_{22} )\\
				&=f (a_{12} x_{22} )+f (b_{12} x_{22} )~\text{(By Lemma \ref{lem3.4} (i))} \\
				&=f (a_{12} ) g (x_{22} )+\alpha (a_{12} ) d (x_{22} )+f (b_{12} ) g (x_{22} )+\alpha (b_{12} ) d (x_{22} )
			\end{aligned}
		\end{equation}
		
		Comparing both sides of \eqref{eq3.9}, we have
		\begin{equation}
			\label{eq3.10}
			\begin{aligned}
				&{ [f (a_{12}+b_{12} )-f (a_{12} )-f (b_{12} ) ] g (x_{22} )=0  }~\text{which implies that,} \\ \\
				&  ( [f (a_{12}+b_{12} )-f (a_{12} )-f (b_{12} ) ]_{11} \\&+ { [f (a_{12}+b_{12} )-f (a_{12} ) }-f (b_{12} ) ]_{12} ) g (x_{22} )=0 \\
				\text{and}~	&  ( [f (a_{12}+b_{12} )-f (a_{12} )-f (b_{12} ) ]_{21} \\&+ { [f (a_{12}+b_{12} )-f (a_{12} ) }-f (b_{12} ) ]_{22} ) g (x_{22} )=0.
			\end{aligned}
		\end{equation}
		
		Let $x_{11} \in R_{11}$. Then
		\begin{equation}
			\label{eq3.11}
			\begin{aligned}
				&f (a_{12}+b_{12} ) g (x_{11} )+\alpha (a_{12}+b_{12} ) d (x_{11} )\\
				&=f ( (a_{12}+b_{12} ) x_{11} ) \\
				&=f(0)=0~\text{(By Lemma \ref{lem3.2})} \\
				&=f (a_{12} x_{11} )+f (b_{12} x_{11} ) \\
				&=f (a_{12} ) g (x_{11} )+\alpha (a_{12} ) d (x_{11} )+f (b_{12} ) g (x_{11} )+\alpha (b_{12} ) d (x_{11} ).
			\end{aligned}
		\end{equation}
		
		Comparing both sides of \eqref{eq3.11},
		\begin{equation}
			\label{eq3.12}
			\begin{aligned}
				&{ [f (a_{12}+b_{12} )-f (a_{12} )-f(b_{12}) ] g (x_{11} )=0 } ~\text{which implies that,} \\ \\
				&  ( [f (a_{12}+b_{12} )-f (a_{12} )-f (b_{12} ) ]_{11} \\&+  { [f (a_{12}+b_{12} )-f (a_{12} ) }-f (b_{12} ) ]_{12} ) g (x_{11} )=0 \\
				\text{and}~	& ( [f (a_{12}+b_{12} )-f(a_{12})-f (b_{12} ) ]_{21}+ \\&{ [f (a_{12}+b_{12} )-f (a_{12} ) }-f (b_{12} ) ]_{22} ) g (x_{11} )=0.
			\end{aligned}
		\end{equation}
		
		By the Assumption \ref{asm1.1} (i), \eqref{eq3.10} and \eqref{eq3.12},
		$$
		\begin{aligned}
			& { [f (a_{12}+b_{12} )-f (a_{12} )-f (b_{12} ) ]_{k 1}=0 } \\
			\& \quad & { [f (a_{12}+b_{12} )-f (a_{12} )-f (b_{12} ) ]_{k 2}=0 . }
		\end{aligned}
		$$
		for $k=1,2$. Hence, we have
		\begin{align*}
			f (a_{12}+b_{12} )=f (a_{12} )+f (b_{12} ).
		\end{align*}	
		
		Similarly, we can prove that
		\begin{align*}
			f (a_{21}+b_{21} )=f (a_{21} )+f (b_{21} ).
		\end{align*}	
	\end{proof}

	\begin{lemma}
		\label{lem3.6}
		Let $a_{ii}, b_{ii} \in R_{ii}$. Then
		$$
		\begin{aligned}
			&{ (i)~ } f (a_{11}+b_{11} )=f (a_{11} )+f (b_{11} ), \\
			& { (ii)~ } f (a_{22}+b_{22} )=f (a_{22} )+f (b_{22} ).
		\end{aligned}
		$$
	\end{lemma}
	
	\begin{proof}
		Let $x_{12} \in R_{12}$. Then
		\begin{equation}
			\label{eq3.13}
			\begin{aligned}
				&f (a_{11}+b_{11} ) g (x_{12} )+\alpha (a_{11}+b_{11} ) d (x_{12} )\\
				&= f ( (a_{11}+b_{11} ) x_{12} ) \\
				&= f (a_{11} x_{12}+b_{11} x_{12} ) \\
				&= f (a_{11} x_{12} )+f (b_{11} x_{12} )~\text{(By Lemma \ref{lem3.5} (i))} \\
				&= f (a_{11} ) g (x_{12} )+\alpha (a_{11} ) d (x_{12} )+f (b_{11} ) g (x_{12} )+\alpha (b_{11} ) d (x_{12} ).
			\end{aligned}
		\end{equation}
		
		Comparing both sides of \eqref{eq3.13},
		\begin{equation}
			\label{eq3.14}
			\begin{aligned}
				&{ [f (a_{11}+b_{11} )-f (a_{11} )-f (b_{11} ) ] g (x_{12} )=0  }~\text{which implies that,} \\ \\
				&  ( [f (a_{11}+b_{11} )-f (a_{11} )-f (b_{11} ) ]_{11} \\&+ { [f (a_{11}+b_{11} )-f (a_{11} ) }-f (b_{11} ) ]_{12} ) g (x_{12} )=0 \\
				\text{and}~	&  ( [f (a_{11}+b_{11} )-f (a_{11} )-f (b_{11} ) ]_{21} \\&+ { [f (a_{11}+b_{11} )-f (a_{11} ) }-f (b_{11} ) ]_{22} ) g (x_{12} )=0.
			\end{aligned}
		\end{equation}

		Similarly, taking $x_{22} \in R_{22}$, we have
		\begin{equation}
			\label{eq3.15}
			\begin{aligned}
				&{ [f (a_{11}+b_{11} )-f (a_{11} )-f (b_{11} ) ] g (x_{22} )=0  }~\text{which implies that,} \\ \\
				&  ( [f (a_{11}+b_{11} )-f (a_{11} )-f (b_{11} ) ]_{11} \\&+ { [f (a_{11}+b_{11} )-f (a_{11} ) }-f (b_{11} ) ]_{12} ) g (x_{22} )=0 \\
			\text{and}~		&  ( [f (a_{11}+b_{11} )-f (a_{11} )-f (b_{11} ) ]_{21} \\&+ { [f (a_{11}+b_{11} )-f (a_{11} ) }-f (b_{11} ) ]_{22} ) g (x_{22} )=0.
			\end{aligned}
		\end{equation}
		
		By the Assumption \ref{asm1.1} (i), \eqref{eq3.14} and \eqref{eq3.15}, we have
		$$
		\begin{aligned}
			&{ [f (a_{11}+b_{11} )-f (a_{11} )-f (b_{11} ) ]_{k1}=0}, \\
			& [f (a_{11}+b_{11} )-f (a_{11} )-f (b_{11} ) ]_{k2}=0 .
		\end{aligned}
		$$
		for $k=1,2$. Hence, we have
		\begin{align*}
			f (a_{11}+b_{11} )=f (a_{11} )+f (b_{11} ).
		\end{align*}	
		
		Similarly, we can prove that
		\begin{align*}
			f (a_{22}+b_{22} )=f (a_{22} )+f (b_{22} ).
		\end{align*}	
	\end{proof}

	\begin{lemma}
		\label{lem3.7}
		Let $a_{ij}\in R_{ij}$. Then
		\begin{align*}
			f (a_{11}+a_{12}+a_{21}+a_{22} )=f (a_{11} )+f (a_{12} )+f (a_{21} )
			+f (a_{22} ).
		\end{align*}
	\end{lemma}
	
	\begin{proof}
		Let $x_{11}\in R_{11}$. Then
		\begin{equation}
			\label{eq3.16}
			\begin{aligned}
				&f (a_{11}+a_{12}+a_{21}+a_{22} ) g (x_{11} )
				+\alpha (a_{11}+a_{12}+a_{21}+a_{22} ) d (x_{11} )\\
				&= f ((a_{11}+a_{12}+a_{21}+a_{22} ) x_{11} ) \\
				&=f (a_{11} x_{11}+a_{21} x_{11} ) \\
				&=f (a_{11} x_{11} )+f (a_{12} x_{11} )+f (a_{21} x_{11} )+f (a_{22} x_{11} ) \\
				&=f (a_{11} ) g (x_{11} )+\alpha (a_{11} ) d (x_{11} )+f (a_{12} ) g (x_{11} )+\alpha (a_{12} ) d (x_{11} ) \\
				&+f (a_{21} ) g (x_{11} )+\alpha (a_{21} ) d (x_{11} )+f (a_{22} ) g (x_{11} ) +\alpha (a_{22} ) d (x_{11} ).
			\end{aligned}
		\end{equation}
		
		Comparing both sides of \eqref{eq3.16},
		\begin{equation}
			\label{eq3.17}
			\begin{aligned}
				& [f (a_{11}+a_{12}+a_{21}+a_{22} )-f (a_{11} )-f (a_{12} ) \\&-f (a_{21}) -f (a_{22} ) ] g (x_{11} )=0  ~\text{which implies that,}\\ \\
				&([f (a_{11}+a_{12}+a_{21}+a_{22} )-f (a_{11} )-f (a_{12} )-f (a_{21}) -f (a_{22} ) ]_{11}\\
				&+[f (a_{11}+a_{12}+a_{21}+a_{22} )-f (a_{11} )-f (a_{12} ) \\&-f (a_{21}) -f (a_{22} ) ]_{12})g (x_{11} )=0\\ \\
				\text{and}~	&([f (a_{11}+a_{12}+a_{21}+a_{22} )-f (a_{11} )-f (a_{12} )-f (a_{21}) -f (a_{22} ) ]_{21}\\
				&+[f (a_{11}+a_{12}+a_{21}+a_{22} )-f (a_{11} )-f (a_{12} ) \\& -f (a_{21}) -f (a_{22} ) ]_{22})g (x_{11} )=0.
			\end{aligned}
		\end{equation}
		
		Similarly, taking $x_{22} \in R_{22}$, we have
		\begin{equation}
			\label{eq3.18}
			\begin{aligned}
				& [f (a_{11}+a_{12}+a_{21}+a_{22} )-f (a_{11} )-f (a_{12} )\\&-f (a_{21}) -f (a_{22} ) ] g (x_{22} )=0 ~\text{which implies that,}\\ \\
				&([f (a_{11}+a_{12}+a_{21}+a_{22} )-f (a_{11} )-f (a_{12} )-f (a_{21}) -f (a_{22} ) ]_{11}\\
				&+[f (a_{11}+a_{12}+a_{21}+a_{22} )-f (a_{11} )-f (a_{12} ) \\&-f (a_{21}) -f (a_{22} ) ]_{12})g (x_{22} )=0\\ \\
				\text{and}~	&([f (a_{11}+a_{12}+a_{21}+a_{22} )-f (a_{11} )-f (a_{12} )-f (a_{21}) -f (a_{22} ) ]_{21}\\
				&+[f (a_{11}+a_{12}+a_{21}+a_{22} )-f (a_{11} )-f (a_{12} ) \\&-f (a_{21}) -f (a_{22} ) ]_{22})g (x_{22} )=0.
			\end{aligned}
		\end{equation}
		
		By the Assumption \ref{asm1.1} (i), \eqref{eq3.17} and \eqref{eq3.18}, we get
		\begin{align*}
			& [f (a_{11}+a_{12}+a_{21}+a_{22} )-f (a_{11} )-f (a_{12} )-f (a_{21}) -f (a_{22} ) ]_{k1}=0,\\
			& [f (a_{11}+a_{12}+a_{21}+a_{22} )-f (a_{11} )-f (a_{12} )-f (a_{21}) -f (a_{22} ) ]_{k2}=0,
		\end{align*}
		for $k=1,2$. Hence, we have
		\begin{align*}
			f (a_{11}+a_{12}+a_{21}+a_{22} )=f (a_{11} )+f (a_{12} )+f (a_{21} )
			+f (a_{22} ).
		\end{align*}
	\end{proof}

	\begin{proof} [Proof of Theorem \ref{thm3.1}] Let $a,b \in R$. Then
		\begin{align*}
			&a=a_{11}+a_{12}+a_{21}+a_{22},\\
			& b= b_{11}+b_{12}+b_{21}+b_{22},
		\end{align*}
		for some $a_{ij}, b_{ij} \in R_{ij}$. Now,
		$$
		\text {  } \begin{aligned}
			f(a+b)=~& f (a_{11}+a_{12}+a_{21}+a_{22}+b_{11}+b_{12}+b_{21}+b_{22} ) \\
			=~& f ( (a_{11}+b_{11} )+ (a_{12}+b_{12} )+ (a_{21}+b_{21} )+ (a_{22}+b_{22} ) ) \\
			=~& f (a_{11}+b_{11} )+f (a_{12}+b_{12} )+f (a_{21}+b_{21} ) +f (a_{22}+b_{22} ) \\
			&~\text{(By Lemma \ref{lem3.7})}\\
			=~& f (a_{11} )+f (b_{11} )+f (a_{12} )+f (b_{12} )+f (a_{21} ) +f (b_{21} ) \\&+f (a_{22} )+f (b_{22} ) ~\text{(By Lemma \ref{lem3.5} and \ref{lem3.6})}\\
			=~& f (a_{11}+a_{12}+a_{21}+a_{22} )+f (b_{11}+b_{12}+b_{21}+b_{22} ) \\
			&~\text{(By Lemma \ref{lem3.7})}\\
			=~& f(a)+f(b).
		\end{aligned}
		$$
		Hence $f$ is additive on $R$. Moreover, $f$ is a generalized skew semi-derivation on $R$.
	\end{proof}

\begin{cor}
	Let $R$ be a ring containing a non-trivial idempotent $e$, $f$ be a multiplicative generalized skew semi-derivation on $R$, $g$ be a map on $R$ with $g(0)=0$ and $f$, $g$ satisfy the following conditions for $i, j, k=1,2~:$ \\
	
	(i) If $(a_{k 1}+b_{k 2}) g(x_{i j})=0$ and $i \leq j$ for all $x_{i j} \in R_{i j}$, then
	$$
	\begin{cases}a_{k 1}=0 & \text { if } i=1, \\ b_{k 2}=0 & \text { if } i=2 .\end{cases}
	$$
	\\
	
	(ii) If $g (x_{i i})(a_{1 j}+b_{2 j})=0$ for all $x_{i i} \in R_{i i}$, then
	$$
	\begin{cases}a_{1 j}=0 & \text { if } i=1, \\ b_{2 j}=0 & \text { if } i=2.\end{cases}
	$$
	\\
	
	(iii) If $f (x_{i i}) (a_{1 j}+b_{2 j})=0$ for all $x_{i i} \in R_{i i}$, then
	$$
	\begin{cases}a_{1 j}=0 & \text { if } i=1, \\ b_{2 j}=0 & \text { if } i=2.\end{cases}
	$$
	Then $f$ is additive.
\end{cor}

\begin{cor}
\label{cor3.9}
	Let $R$ be a prime ring with identity $1\neq 0$ and a non-trivial idempotent $e$. If $f$ is a multiplicative generalized skew semi-derivation on $R$ where the associated map $g$ is an endomorphism on $R$ such that $g(e)=e$, then $f$ is additive.
\end{cor}

\begin{cor}
Let $R$ be a ring as in Corollary \ref{cor3.9}. If $f$ is a multiplicative generalized skew  semi-derivation on $R$, where the associated map $g$ is invariant on $R_{i j}$, i.e., $g (R_{i j}) \subseteq R_{i j}$ where $i, j=1,2$ and satisfies conditions $(iii)$ and $(iv)$ of Assumption \ref{asm1.1}, then $f$ is additive.
\end{cor}

\begin{theorem}
\label{thm3.11}
Let $R$ be a ring with a non-trivial idempotent $e$, $g$ be a map on $R$ vanishing at zero and satisfies the following conditions for $i, j=1,2~:$\\
	
	(i) If $ (a_{j 1}+b_{j 2}) g (x_{i i})=0$ for all $x_{i i} \in R_{i i}$, then $a_{j 1}+b_{j 2}=0$.\\
	
	(ii) If $g (x_{i i})(a_{11}+a_{21})=0$ for all $x_{i i} \in R_{i i}$, then $a_{11}+a_{21}=0$.\\
	
	If $f$ is a multiplicative generalized skew semi-derivation on $R$, then $f$ is additive.
\end{theorem}

	\begin{cor}
		Let $R$ be a ring containing a non-trivial idempotent $e$, $f$ be a multiplicative generalized skew semi-derivation on $R$, $g$ be a map on $R$ such that $g(0)=0$ and satisfies three of the following conditions for $i=1,2~:$ \\
		
		(i) If $ (a_{11}+b_{12}) g (x_{ii})=0$ for all $x_{ii} \in R_{ii}$, then $a_{11}+b_{12}=0$,\\
		
		(ii) If $ (a_{21}+b_{22}) g (x_{ii})=0$ for all $x_{ii} \in R_{ii}$, then $a_{21}+b_{22}=0$,\\
		
		(iii) If $g (x_{ii}) (a_{11}+b_{21})=0$ for all $x_{ii} \in R_{1 1}$, then $a_{11}+b_{21}=0$,\\
		
		(iv) If $g(x_{ii}) (a_{12}+b_{22})=0$ for all $x_{ii} \in R_{ii}$, then $a_{12}+b_{22}=0$.\\
		
		Then $f$ is additive.	
	\end{cor}
\section*{Acknowledgement}
The first author is thankful to the University Grants Commission (UGC), Govt. of India for financial supports under File. No. 16-9(June 2019)/2019(NET/CSIR), UGC Ref. No. 1256/(CSIR-UGC NET JUNE 2019) dated 16/12/2019 and all authors are thankful to the Indian Institute of Technology Patna for providing the research facilities.
\section*{Declarations}
\textbf{Data Availability Statement}: The authors declare that [the/all other] data supporting the findings of this study are available within the article. \\
\textbf{Competing interests}: The authors declare that there is no conflict of interest regarding the publication of this manuscript.\\

\end{document}